\documentclass{amsart}
\usepackage{latexsym}
\usepackage{amssymb}
\usepackage{amsmath}
\usepackage[centredisplay,PostScript=dvips]{diagrams}

\def\gt{\mathfrak{T}}
\newtheorem{df}{Definition}[section]

\newtheorem{lem}[df]{Lemma}
\newtheorem{coro}[df]{Corollary}

\begin{document}

\title[A formula for the First Eigenvalue of the Dirac Operator...]
{A formula for the First Eigenvalue of the Dirac Operator on
Compact Spin Symmetric Spaces}%
\author{Jean-Louis Milhorat}
\address{Laboratoire Jean Leray\\
 UMR CNRS 6629\\
 D\'epartement de Math\'ematiques\\
 Universit\'e de Nantes\\
 2, rue de la Houssini\`ere\\
 BP 92208\\
  F-44322 NANTES CEDEX 03}
\email{milhorat@math.univ-nantes.fr}


\date{}
\begin{abstract}
Let $G/K$ be a simply connected spin compact inner irreducible
symmetric space, endowed with the metric induced by the Killing
form of $G$ sign-changed. We give a formula for the square of the
first eigenvalue of the Dirac operator in terms of a root system
of $G$. As an example of application, we give the list of the
first eigenvalues for the spin compact irreducible symmetric
spaces endowed with a quaternion-K\"ahler structure.
\end{abstract}

\maketitle
\section{Introduction}
Let $G/K$ be a compact, simply-connected, $n$-dimensional
irreducible symmetric space with $G$ compact and simply-connected,
endowed with the metric induced by the Killing form of $G$
sign-changed. Assume that $G$ and $K$ have same rank and that
$G/K$ has a spin structure. In a previous paper, cf. \cite{Mil04},
we proved that the first eigenvalue $\lambda$ of the Dirac
operator verifies
\begin{equation}
\label{formula1}\lambda^{2}=2\,\min_{1\leq k\leq
p}\|\beta_{k}\|^{2}+n/8\,,
\end{equation}
where $\beta_{k}$, $k=1,\ldots, p$, are the $K$-dominant weights
occurring in the decomposition into irreducible components of the
spin representation under the action of $K$, and where $\|\cdot\|$
is the norm associated to the scalar product induced by the
Killing form of $G$.

\noindent The proof was based on a lemma of R. Parthasarathy in
\cite{Par}, which allows to express the result in the following
way.

\noindent Let $T$ be a fixed common maximal torus of $G$ and $K$.
Let $\Phi$ be
  the set of non-zero roots of $G$ with respect to $T$.
  Let $\Phi_{G}^{+}$  be the set of positive roots of $G$,
$\Phi_{K}^{+}$
 be the set of positive roots of $K$, with
respect to a fixed lexicographic ordering in $\Phi$. Let
$\delta_{G}$, (resp. $\delta_{K}$) be the half-sum of the positive
roots of $G$, (resp. $K$). Then the square of the first eigenvalue
of the Dirac operator is given by \begin{equation} \lambda^{2}=
2\,\min_{w\in W}\|w\cdot \delta_{G}-\delta_{K}\|^{2}+n/8\,,
\end{equation} where $W$ is the subset of the Weyl group $W_{G}$ defined
by
\begin{equation}\label{def.w} W:=\{w\in W_{G}\;;\;
w\cdot \Phi_{G}^{+}\supset \Phi_{K}^{+}\}\,.
\end{equation}
In order to avoid the determination of the subset $W$ for
applications, we prove in the following that the square of the
first eigenvalue of the Dirac operator is indeed given by
\begin{equation}\label{formula2} \boxed{\lambda^{2}=
2\,\min_{w\in W_{G}}\|w\cdot \delta_{G}-\delta_{K}\|^{2}+n/8\,.}
\end{equation}

\noindent We then give a different expression to use the formula
for explicit computations. We obtain
\begin{equation}\label{formula3}
\boxed{\lambda^{2}=2\,\|\delta_{G}-\delta_{K}\|^{2}+4\,\sum_{\theta\in\Lambda}
<\theta,\delta_{K}>+n/8\,,}\end{equation} where $\Lambda$ is the
set $$\Lambda:=\{\theta\in \Phi_{G}^{+}\,;\,
<\theta,\delta_{K}>\;\;<0\}\,.$$

\noindent As an example of application of the above formula, we
obtain the list of the first eigenvalues of the Dirac operator for
the spin compact irreducible symmetric spaces endowed with a
quaternion-K\"ahler structure. By definition, a Riemannian
manifold has a quaternion-K\"ahler structure if its holonomy group
is contained in the group $\mathrm{Sp}_{m}\mathrm{Sp}_{1}$. In
\cite{Wol65}, J. Wolf gave the following classification of compact
quaternion-K\"ahler symmetric spaces:
 \vspace{5mm}
\begin{center}
\begin{tabular} {|c|c|c|c|c|} \hline
 $G$&$K$&$ G/K$& $\dim\, G/K$&Spin structure  \\ & & & &(cf.
\cite{CG})\\ \hline $\mathrm{Sp}_{m+1}$
&$\mathrm{Sp}_{m}\times\mathrm{Sp}_{1}$& Quaternionic & $4m\,
(m\geq 1)$&Yes (unique) \\ &&projective&&\\ &&space
$\mathbb{H}P^{m}$&&\\ \hline
$\mathrm{SU}_{m+2}$&$S(\mathrm{U}_{m}\times\mathrm{U}_{2})$ &
Grassmannian &$4m\, (m\geq 1)$&iff $m$ even
\\
&&$\mathrm{Gr}_{2}(\mathbb{C}^{m+2})$&&unique in that case\\
\hline$\mathrm{Spin}_{m+4}$& $\mathrm{Spin}_{m}\mathrm{Spin}_{4}$&
Grassmannian &$4m\, (m\geq 3)$ &iff $m$ even,
\\
&&$\widetilde{\mathrm{Gr}}_{4}(\mathbb{R}^{m+4})$&&unique in that
case\\ \hline $\mathrm{G}_{2}$&$\mathrm{SO}_{4}$&&$8$&Yes
(unique)\\&&&&\\ \hline $\mathrm{F}_{4}$ &
$\mathrm{Sp}_{3}\mathrm{SU}_{2}$&&$28$&No\\&&&&\\ \hline $
\mathrm{E}_{6}$& $\mathrm{SU}_{6}\mathrm{SU}_{2}$&&$40$&Yes
(unique)\\&&&&\\ \hline $\mathrm{E}_{7}$&
$\mathrm{Spin}_{12}\mathrm{SU}_{2}$&&$64$&Yes (unique)\\ &&&&\\
\hline $\mathrm{E}_{8}$ & $\mathrm{E}_{7}\mathrm{SU}_{2}$&&$112$&
Yes (unique)\\&&&&\\ \hline
\end{tabular}\end{center}

\vspace{5mm}

 \noindent Note furthermore that all the symmetric spaces in that
 list are ``{inner}''.

\noindent Endowing each symmetric space with the metric induced by
the Killing form of $G$ sign-changed, we obtain the following
table

\vfill \newpage

\begin{tabular}{|c|c|}\hline
&\\ $G/K$ & Square of the first eigenvalue of D\\&\\ \hline &\\
$\mathbb{H}P^n=\mathrm{Sp}_{m+1}/(\mathrm{Sp}_{m}\times\mathrm{Sp}_{1})$&
$\displaystyle{\frac{m+3}{m+2}\,\frac{m}{2}=\frac{m+3}{m+2}\,\frac{\mathrm{Scal}}{4}}$\\
&\\\hline&\\ $\mathrm{Gr}_{2}(
\mathbb{C}^{m+2})=\mathrm{SU}_{m+2}/S(\mathrm{U}_{m}\times\mathrm{U}_{2})$&
$\displaystyle{\frac{m+4}{m+2}\,\frac{m}{2}=\frac{m+4}{m+2}\,\frac{\mathrm{Scal}}{4}}$\\
($m$ even)&\\ \hline &\\ $\widetilde{\mathrm{Gr}}_{4}(
\mathbb{R}^{m+4})=\mathrm{Spin}_{m+4}/\mathrm{Spin}_{m}\mathrm{Spin}_{4}$&
$\displaystyle{\frac{m^2+6m-4}{m(m+2)}\,\frac{m}{2}=\frac{m^2+6m-4}{m(m+2)}
\,\frac{\mathrm{Scal}}{4}}$\\ ($m$ even) &\\ \hline &\\
$\mathrm{G}_{2}/\mathrm{SO}_{4}$&$\displaystyle{\frac{3}{2}=
\frac{3}{2}\,\frac{\mathrm{Scal}}{4}}$\\ &\\ \hline &\\ $
\mathrm{E}_{6}/(\mathrm{SU}_{6}\mathrm{SU}_{2})$ &
$\displaystyle{\frac{41}{6}=\frac{41}{30}\,\frac{\mathrm{Scal}}{4}}$\\
&\\ \hline &\\ $\mathrm{E}_{7}/(\mathrm{Spin}_{12}\mathrm{SU}_{2})
$ &
$\displaystyle{\frac{95}{9}=\frac{95}{72}\,\frac{\mathrm{Scal}}{4}}$\\
&\\ \hline &\\ $\mathrm{E}_{8}/(\mathrm{E}_{7}\mathrm{SU}_{2})$ &
$\displaystyle{\frac{269}{15}=\frac{269}{210}\,\frac{\mathrm{Scal}}{4}}$\\&\\
\hline
\end{tabular}
\vspace{5mm}

\centerline{TABLE I}

\vspace{5mm}

\noindent The result was already known for quaternionic projective
spaces $\mathbb{H}P^n$, \cite{Mil92}, for the Grassmannians
$\mathrm{Gr}_{2}( \mathbb{C}^{m+2})$, \cite{Mil98}, and for the
symmetric space $\mathrm{G}_{2}/\mathrm{SO}_{4}$, \cite{See99}. Up
to our knowledge, the other results are new.

\section{Proof of formula (\ref{formula2})}
With the notations of the introduction, and since the scalar
product is $W_{G}$-invariant, one has for any $w\in W_{G}$
\begin{equation}\label{eq1}\|w\cdot \delta_{G}-\delta_{K}\|^{2}=\|\delta_{G}\|^{2}+
\|\delta_{K}\|^{2}-2\,<w\cdot
\delta_{G},\delta_{K}>\,,\end{equation} hence
\begin{align*}
\min_{w\in W}\|w\cdot \delta_{G}-\delta_{K}\|^{2}&=
\|\delta_{G}\|^{2}+ \|\delta_{K}\|^{2}-2\,\max_{w\in W}\,<w\cdot
\delta_{G},\delta_{K}>,,\\ \intertext{and} \min_{w\in
W_{G}}\|w\cdot \delta_{G}-\delta_{K}\|^{2}&= \|\delta_{G}\|^{2}+
\|\delta_{K}\|^{2}-2\,\max_{w\in W_{G}}\,<w\cdot
\delta_{G},\delta_{K}>,.
\end{align*}
So we have to prove that \begin{equation}\max_{w\in W}\,<w\cdot
\delta_{G},\delta_{K}>=\max_{w\in W_{G}}\,<w\cdot
\delta_{G},\delta_{K}>\,.\end{equation} Let
\begin{equation}\Pi_{G}:=\{\theta_{1},\ldots,\theta_{r}\}\subset
\Phi_{G}^{+}\,,\end{equation} be the set of $G$-simple roots and
let
\begin{equation}\Pi_{K}:=\{\theta_{1}',\ldots,\theta_{l}'\}\subset
\Phi_{K}^{+}\,,\end{equation} be the set of $K$-simple roots.

\noindent Let $w_{0}\in W_{G}$ such that
\begin{equation}\label{w0}<w_{0}\cdot \delta_{G},\delta_{K}>=\max_{w\in
W_{G}}\,<w\cdot \delta_{G},\delta_{K}>\,. \end{equation} Suppose
that $w_{0}\notin W$. Then we claim that there exists a $K$-simple
root $\theta_{i}'$ such that $w_{0}^{-1}\cdot \theta_{i}' \notin
\Phi_{G}^{+}$. Otherwise, if for any $K$-simple root
$\theta_{i}'$, $w_{0}^{-1}\cdot \theta_{i}' \in \Phi_{G}^{+}$,
then since any $K$-positive root is a linear combination with
non-negative coefficients of $K$-simple roots, we would have
$\forall \theta' \in \Phi_{K}^{+}$, $w_{0}^{-1}\cdot\theta' \in
\Phi_{G}^{+}$, contradicting the assumption made on $w_{0}$.

\noindent Now let $\sigma_{i}'$ be the reflection across the
hyperplane ${\theta_{i}'}^{\bot}$. Since
$\sigma_{i}'\cdot\delta_{K}=\delta_{K}-\theta_{i}'$, (cf. for
instance Corollary of Lemma~B, \S10~.3 in \cite{Hum}), one gets by
the $W_{G}$-invariance of the scalar product
\begin{align*}
<\sigma_{i}'w_{0}\cdot \delta_{G},\delta_{K}>&=<w_{0}\cdot
\delta_{G}, \sigma_{i}'\cdot\delta_{K}>=<w_{0}\cdot
\delta_{G},\delta_{K}-\theta_{i}'>\\ &= <w_{0}\cdot
\delta_{G},\delta_{K}>-<\delta_{G},w_{0}^{-1}\cdot\theta_{i}'>\,.
\end{align*}
But since $w_{0}^{-1}\cdot\theta_{i}'$ is a negative root of $G$,
one has $$w_{0}^{-1}\cdot\theta_{i}'=-\sum k_{j}\,
\theta_{j}\,,\quad k_{j}\in \mathbb{N}\,.$$ Since for any
$G$-simple root $\theta_{j}$,
$\sigma_{j}\cdot\delta_{G}=\delta_{G}-\theta_{j}$, where
$\sigma_{j}$ is the reflection across the hyperplane
${\theta_{j}}^{\bot}$, one has $<\theta_{j},\delta_{G}>=2\,
<\theta_{j},\theta_{j}>\; >0$, so $$
-<\delta_{G},w_{0}^{-1}\cdot\theta_{i}'>=\sum
k_{j}\,<\delta_{G},\theta_{j}>\; >0\,,$$ hence $$
<\sigma_{i}'w_{0}\cdot \delta_{G},\delta_{K}>\; >\; <w_{0}\cdot
\delta_{G},\delta_{K}>\,,$$ but that is in contradiction with the
definition (\ref{w0}) of $w_{0}$, hence $w_{0}\in W$ and
$$\max_{w\in W_{G}}\,<w\cdot \delta_{G},\delta_{K}>=<w_{0}\cdot
\delta_{G},\delta_{K}> \leq \max_{w\in W}\,<w\cdot
\delta_{G},\delta_{K}>\leq \max_{w\in W_{G}}\,<w\cdot
\delta_{G},\delta_{K}>\,,$$ hence the result.

\section{Proof of formula (\ref{formula3})}
In order to obtain the formula we will use the following result
\begin{lem}
For any element $w$ of the Weyl group $W_{G}$
$$w\cdot\delta_{G}=\delta_{G}-\sum_{\theta\in\Phi_{G}^{+}}k_{\theta}\,\theta
\,,\quad k_{\theta}=0\text{ or }\,1 \,.$$
\end{lem}
\begin{proof} Let $w\in W_{G}$. With the same notations as in the
above proof, we write $w$ in reduced form
\begin{equation}\label{red.exp}w=\sigma_{i_{1}}\cdots\sigma_{i_{k}}\,,
\end{equation} where $\sigma_{i}$ is the reflection across the hyperplane
$\theta_{i}^{\bot}$, $\theta_{i}\in \Pi_{G}$, and $k$ is minimal.

\noindent Since
$\sigma_{i_{k}}\cdot\delta_{G}=\delta_{G}-\theta_{i_{k}}$, one has
$$w\cdot\delta_{G}=
\sigma_{i_{1}}\cdots\sigma_{i_{k-1}}(\sigma_{i_{k}}\cdot\delta_{G})=
\sigma_{i_{1}}\cdots\sigma_{i_{k-1}}(\delta_{G}) -
\sigma_{i_{1}}\cdots\sigma_{i_{k-1}}(\theta_{i_{k}})\,.$$ Now,
since the expression of $w$ is reduced, $w(\theta_{i_{k}})$ is a
negative root, cf. for instance corollary of Lemma~C, \S~10.3 in
\cite{Hum}. But
$w(\theta_{i_{k}})=-\sigma_{i_{1}}\cdots\sigma_{i_{k-1}}(\theta_{i_{k}})$,
hence $\sigma_{i_{1}}\cdots\sigma_{i_{k-1}}(\theta_{i_{k}})$ is a
positive root.

\noindent Now the element $\sigma_{i_{1}}\cdots\sigma_{i_{k-1}}\in
W_{G}$
 is written in reduced form, otherwise the
expression (\ref{red.exp}) of $w$ would not be reduced. Hence we
may conclude as above that
$$\sigma_{i_{1}}\cdots\sigma_{i_{k-1}}(\delta_{G})=
\sigma_{i_{1}}\cdots\sigma_{i_{k-2}}(\delta_{G})-
\sigma_{i_{1}}\cdots\sigma_{i_{k-2}}(\theta_{i_{k-1}})\,,$$ where
$ \sigma_{i_{1}}\cdots\sigma_{i_{k-2}}(\theta_{i_{k-1}})$ is a
positive root.

\noindent Proceeding inductively we get
$$w\cdot\delta_{G}=\delta_{G}-\sum_{\theta\in\Phi_{G}^{+}}k_{\theta}\,\theta
\,,\quad k_{\theta}\in \mathbb{N} \,.$$ In order to conclude, we
have to prove that if a $G$-positive root $\theta$ appears in the
above sum, then it appears only once.

\noindent Suppose that a $G$-positive root appears at least twice
in the above sum, then there exist two integers $p$ and $q$,
$1\leq p<q\leq k-1$ such that
$$\sigma_{i_{1}}\cdots\sigma_{i_{p}}(\theta_{i_{p+1}})
=\sigma_{i_{1}}\cdots\sigma_{i_{q}}(\theta_{i_{q+1}})\,.$$
applying $\sigma_{i_{p+1}}\sigma_{i_{p}}\cdots \sigma_{i_{1}}$ to
the two members of the above equation, we get
$$\begin{cases}-\theta_{i_{p+1}}= \sigma_{i_{p+2}}\cdots
\sigma_{i_{q}}(\theta_{i_{q+1}})\,, &\text{if $p+1<q$,}\\
-\theta_{i_{q}}=\theta_{i_{q+1}}\,, &\text{if $p+1=q$.}
\end{cases}$$
So we get a contradiction, even in the first case, since
$\sigma_{i_{p+2}}\cdots \sigma_{i_{q}}\sigma_{i_{q+1}}\in W_{G}$
is expressed in reduced form (otherwise the expression
(\ref{red.exp}) of $w$ would not be reduced), hence
$\sigma_{i_{p+2}}\cdots \sigma_{i_{q}}(\theta_{i_{q+1}})$ is a
positive root.
\end{proof}

\noindent From the above result we deduce
\begin{lem}\label{lemLambda}
Let $\Lambda$ be the set \begin{equation}\Lambda:=\{\theta\in
\Phi_{G}^{+}\,;\, <\theta,\delta_{K}>\;\;<0\}\,.\end{equation} One
has $$ \max_{w\in W_{G}}\,<w\cdot
\delta_{G},\delta_{K}>=<\delta_{G},\delta_{K}>-\sum_{\theta\in\Lambda}
<\theta,\delta_{K}>\,, $$ (setting $\sum_{\theta\in\Lambda}
<\theta,\delta_{K}>=0$, if $\Lambda=\emptyset$).
\end{lem}
\begin{proof} Suppose $\Lambda\neq\emptyset$. We first prove that there
exists $w_{0}\in W_{G}$ such that
$$w_{0}\cdot\delta_{G}=\delta_{G}-\sum_{\theta\in\Lambda}
\theta\,.$$ Let
$$\Phi_{n}^{+}:=\Phi_{G}^{+}\backslash\Phi_{K}^{+}\,.$$ We first
remark that any root in $\Lambda$ belongs to $\Phi_{n}^{+}$.
Otherwise, if there exists $\theta\in \Lambda \cap \Phi_{K}^{+}$,
then since $\theta$ is a combination with non-negative
coefficients of simple $K$-roots, and since
$<\delta_{K},\theta_{i}'>\,>\, 0$, for any $K$-simple root
$\theta_{i}'$, we would have $<\delta_{K},\theta>\geq 0$,
contradicting the fact that $\theta\in \Lambda$.

\noindent Now, consider $$\delta_{n}:=\frac{1}{2}\sum_{\theta \in
\Phi_{n}^{+}}\theta= \delta_{G}-\delta_{K}\,.$$ Then
$$\delta_{G}-\sum_{\theta\in\Lambda}
\theta=\delta_{K}+\left(\delta_{n}-\sum_{\theta\in\Lambda}\theta\right)\,.$$
But, $$\beta:=\delta_{n}-\sum_{\theta\in\Lambda}\theta\,,$$ is a
weight of the decomposition of the spin representation under the
action of $K$, cf. \S~2 in \cite{Par}: the weights are just the
elements of the form $\delta_{n}-\sum_{\theta\in \Upsilon}\theta$,
where $\Upsilon$ is a subset of $\Phi_{n}^{+}$.

\noindent In fact $\beta$ is the highest weight of an irreducible
component in the decomposition, otherwise we would have $$\beta
+\alpha=\delta_{n}-\sum_{\theta\in \Upsilon}\theta\,,$$ where
$\alpha$ is a $K$-positive root and $\Upsilon$ is a subset of
$\Phi_{n}^{+}$.

\noindent Hence setting $\Lambda':=\Lambda\backslash\Upsilon$ and
$\Upsilon':=\Upsilon\backslash\Lambda$, we would have
$$-\sum_{\theta\in\Lambda'}\theta+\alpha=-\sum_{\theta\in\Upsilon'}\theta\,.$$
But since $\Lambda'\subset\Lambda$ and $\alpha$ is a $K$-positive
root
$$<-\sum_{\theta\in\Lambda'}\theta+\alpha,\delta_{K}>\;>\;0\,,$$
whereas since $\Upsilon'\subset \Phi_{n}^{+}\backslash\Lambda$
$$<-\sum_{\theta\in\Upsilon'}\theta,\delta_{K}>\leq 0\,,$$ hence a
contradiction.

\noindent Now by the result of lemma~2.2 in \cite{Par}, any
highest weight in the decomposition of the spin representation has
the form $$ w\cdot\delta_{G}-\delta_{K}\,,$$ where $w$ belongs to
the subset $W$ of $W_{G}$ defined in (\ref{def.w}). Hence there
exists a $w_{0}\in W$ such that
$$\beta=w_{0}\cdot\delta_{G}-\delta_{K}\,,$$ hence
$$\delta_{G}-\sum_{\theta\in\Lambda}
\theta=\delta_{K}+\beta=w_{0}\cdot\delta_{G}\,,$$ hence the
result.

\noindent Now let $w$ be any element in $W_{G}$. By the above
lemma,
\begin{align*}w\cdot\delta_{G}&=\delta_{G}-\sum_{\theta\in\Phi_{G}^{+}}k_{\theta}\,\theta
\,,&& k_{\theta}=0\text{ or }\,1 \,,\\ &=
\delta_{G}-\sum_{\theta\in\Lambda}k_{\theta}\,\theta-
\sum_{\theta\in\Phi_{G}^{+}\backslash\Lambda}k_{\theta}\,\theta\,.&&
\end{align*}
Hence by the definition of $\Lambda$
$$<w\cdot\delta_{G},\delta_{K}>\leq
<\delta_{G}-\sum_{\theta\in\Lambda}k_{\theta}\,\theta,\delta_{K}>\leq
<\delta_{G}-\sum_{\theta\in\Lambda}\theta,\delta_{K}>\,.$$ Thus
\begin{multline}\max_{w\in W_{G}}\,<w\cdot \delta_{G},\delta_{K}>\leq
<\delta_{G},\delta_{K}>-\sum_{\theta\in\Lambda}
<\theta,\delta_{K}>=<w_{0}\cdot \delta_{G},\delta_{K}> \\ \leq
\max_{w\in W_{G}}\,<w\cdot \delta_{G},\delta_{K}>
\,,\nonumber\end{multline} hence the result.
\end{proof}
 \noindent Now going back to formula~(\ref{formula2}), we get
immediately from (\ref{eq1})
\begin{coro} The first eigenvalue $\lambda$ of the Dirac operator
verifies
$$\lambda^{2}=2\,\|\delta_{G}-\delta_{K}\|^{2}+4\,\sum_{\theta\in\Lambda}
<\theta,\delta_{K}>+n/8\,.$$
\end{coro}

\section{Proof of the results of Table I}

In the following, we note for any integer $n\geq 1$,
$(e_{1},\ldots,e_{n})$, the standard basis of $\mathbb{K}^{n}$,
$\mathbb{K}=\mathbb{R}$, $\mathbb{C}$ or $\mathbb{H}$. The space
of $(n,n)$ matrices with coefficients in $\mathbb{K}$ is denoted
by $\mathrm{M}_{n}(\mathbb{K})$.

\subsection{Quaternionic projective spaces $\boldsymbol{\mathbb{H}P^{n}}$}
 Here $G=\mathrm{Sp}_{m+1}$ and
$K=\mathrm{Sp}_{m}\times\mathrm{Sp}_{1}$. The decomposition of the
spin representation into irreducible components under the action
of $K$ is given in \cite{Mil92}, so we may conclude with
formula~(\ref{formula1}). However the result may be also simply
concluded with formula~(\ref{formula3}).

\noindent The space $\mathbb{H}^{n+1}$ is viewed as a right vector
space on $\mathbb{H}$ in such a way that $G$ may be identified
with the group $$ \left\{ A\in
\mathrm{M}_{m+1}(\mathbb{H})\,;\,{}^{t}AA=I_{m+1}\right\}\,,$$
acting on the left on $\mathbb{H}^{n+1}$ in the usual way. The
group $K$ is identified with the subgroup of $G$ defined by $$
\left\{ A\in \mathrm{M}_{m+1}(\mathbb{H})\,;\,A=\begin{pmatrix}
B&0\\0&q
\end{pmatrix}\,,\,{}^{t}BB=I_{m}\,,\, q\in
\mathrm{Sp}_{1}\right\}\,.$$ Let $T$ be the common torus of $G$
and $K$ $$ T:=\left\{\begin{pmatrix}
\mathrm{e}^{\mathbf{i}\beta_{1}}&&\\ &\ddots&\\ &&
\mathrm{e}^{\mathbf{i}\beta_{m+1}}\end{pmatrix}\;,\;  \beta_{1}
,\ldots, \beta_{m+1}\, \in \mathbb{R} \right\}\; ,$$ where
$$\forall \beta \in \mathbb{R}\,,\quad
\mathrm{e}^{\mathbf{i}\beta}:=\cos(\beta)+\sin(\beta)\,\mathbf{i}\,,$$
$(1,\mathbf{i},\mathbf{j},\mathbf{k})$ being the standard basis of
$\mathbb{H}$.

\noindent The Lie algebra of $T$ is $$\gt=\left\{
\begin{pmatrix} \mathbf{i}\beta_{1}&&\\ &\ddots&\\ &&
\mathbf{i}\beta_{m+1}\end{pmatrix}\; ;\; \beta_{1}, \beta_{2},
,\ldots, \beta_{m+1}\, \in \mathbb{R} \right\}\; .$$ We denote by
$(x_1, \ldots , x_{m+1})$ the  basis of  ${\gt}^{*}$ given by $$
x_k\cdot\begin{pmatrix} \mathbf{i}\beta_{1}&&\\ &\ddots&\\ &&
\mathbf{i}\beta_{m+1}\end{pmatrix}= \beta_{k}\, .$$ A vector
$\mu\in i\, \gt^{*}$ such that $\mu=\sum_{k=1}^{m+1} \mu_{k}\,
\widehat{x}_{k}$, in the basis  $(\widehat {x}_k\equiv i\,
x_k)_{k=1,\ldots ,m+1}$, is denoted by $$ \mu=
(\mu_{1},\mu_{2},\ldots ,\mu_{m+1})\, .$$ The restriction to $\gt$
of the Killing form $B$ of $G$ is given by $$\forall
X\in\gt\,,\,\forall Y \in\gt\,,\quad B(X,Y)= 4\, (m+2)\,\Re\,
\big(\,\mathrm{tr} (X\,Y)\,\big)\,.$$ It is easy to verify that
the scalar product on $i \, \gt^{*}$ induced by the Killing form
sign changed is given by \begin{equation}\begin{split}
\forall\mu=(\mu_1,\ldots ,\mu_{m+1})\in i \, \gt^{*}\,&,\,\forall
\mu'=(\mu'_1,\ldots , \mu'_{m+1})\, \in i\, \gt^*\,,\\ <\mu
,\mu'>&=\frac{1}{4(m+2)}\,\sum_{k=1}^{m+1} \mu_k\, \mu'_k
\,.\end{split}\end{equation} Now, considering the decomposition of
the complexified Lie algebra of $G$ under the action of $T$, it is
easy to verify that $T$ is a common maximal torus of $G$ and $K$,
and that the respective roots are given by
\begin{align*}
&\begin{cases}\pm(\widehat{x}_{i}+\widehat{x}_{j})\,,\\
\pm(\widehat{x}_{i}-\widehat{x}_{i})\,,\end{cases}\; 1\leq i<
j\leq m+1\,,&&\pm 2\,\widehat{x}_{i}\,,\; 1\leq i\leq m+1
&&\text{for}\; G\,,\\
&\begin{cases}\pm(\widehat{x}_{i}+\widehat{x}_{j})\,,\\
\pm(\widehat{x}_{i}-\widehat{x}_{j})\,,\end{cases}\; 1\leq i<
j\leq m\,,&&\pm 2\,\widehat{x}_{i}\,,\; 1\leq i\leq m+1
&&\text{for}\; K\,.
\end{align*}
We consider as sets of positive roots
\begin{align*}\Phi_{G}^{+}&=\left\{\begin{cases}\widehat{x}_{i}+\widehat{x}_{j}\,,\\
\widehat{x}_{i}-\widehat{x}_{j}\,,\end{cases} 1\leq i\leq j\leq
m+1\,;\; 2\,\widehat{x}_{i}\,,\; 1\leq i \leq m+1\right\}\,,\\
\intertext{and}
\Phi_{K}^{+}&=\left\{\begin{cases}\widehat{x}_{i}+\widehat{x}_{j}\,,\\
\widehat{x}_{i}-\widehat{x}_{j}\,,\end{cases}1\leq i\leq j\leq
m\,;\; 2\,\widehat{x}_{i}\,,\; 1\leq i\leq m+1\right\}\,.
\end{align*}
Then
\begin{align*}
\delta_{G}&=\sum_{k=1}^{m+1}
(m+2-k)\,\widehat{x}_{k}=(m+1,m,\ldots,2,1)\,,\\ \intertext{and}
\delta_{K}&=\sum_{k=1}^{m}
(m+1-k)\,\widehat{x}_{k}+\widehat{x}_{m+1}=(m,m-1,\ldots,1,1)\,.
\end{align*}
Hence
$$\delta_{G}-\delta_{K}=\sum_{k=1}^{m}\widehat{x}_{k}=(1,1,\ldots,1,0)\,,$$
so $$\|\delta_{G}-\delta_{K}\|^{2}=\frac{m}{4(m+2)}\,.$$ On the
other hand, it is easy to verify that the set
$$\Lambda:=\{\theta\in \Phi_{G}^{+}\,;\,
<\theta,\delta_{K}>\;\;<0\}\,,$$ is empty, hence by
formula~(\ref{formula3}), the square of the first eigenvalue
$\lambda$ of the Dirac operator is given by $$\lambda^{2}=
\frac{m}{2(m+2)}+\frac{m}{2}=\frac{m+3}{m+2}\,\frac{m}{2}\,.$$

\subsection{Grassmannians $\boldsymbol{\mathrm{Gr}_{2}(\mathbb{C}^{m+2})}$,
$\boldsymbol{m}$ even $\boldsymbol{\geq 2}$} Here
$G=\mathrm{SU}_{m+2}$ and $K$ is the subgroup
$S(\mathrm{U}_{m}\times\mathrm{U}_{2})$ defined below. Here again,
the decomposition into irreducible components of the spin
representation under the action of $K$ is known, \cite{Mil98},
hence the result may be obtained from formula~(\ref{formula1}).
However the result may be also simply concluded with
formula~(\ref{formula3}).

\noindent The group $G$ is identified with $$\left\{ A\in
\mathrm{M}_{m+2}(\mathbb{C})\,;\, {}^{t}AA=I_{m+2}\text{ and }
\det A=1\right\}\,.$$ The group $K$ is the group
$$S(\mathrm{U}_{m}\times\mathrm{U}_{2})=\left\{ A\in
\mathrm{M}_{m+2}(\mathbb{C})\,;\,A=\begin{pmatrix} B&0\\0&C
\end{pmatrix}\,,\,B\in\mathrm{U}_{m}\,,\,C\in\mathrm{U}_{2}\,;\,
\det A=1 \right\}\,.$$ Let $T$ be the common torus of $G$ and $K$
$$ T:=\left\{\begin{pmatrix} \mathrm{e}^{i\beta_{1}}&&\\
&\ddots&\\ && \mathrm{e}^{i\beta_{m+2}}\end{pmatrix}\;,\;
\beta_{1} ,\ldots, \beta_{m+2}\, \in
\mathbb{R}\,,\,\sum_{k=1}^{m+2}\beta_{k}=0 \right\}\; .$$

\noindent The Lie algebra of $T$ is $$\gt=\left\{
\begin{pmatrix} i\beta_{1}&&\\ &\ddots&\\ &&
i\beta_{m+2}\end{pmatrix}\; ;\; \beta_{1}, \beta_{2}, ,\ldots,
\beta_{m+2}\, \in \mathbb{R}\,,\,\sum_{k=1}^{m+2}\beta_{k}=0
\right\}\; .$$ We denote by $(x_1, \ldots , x_{m+1})$ the  basis
of  ${\gt}^{*}$ given by $$ x_k\cdot\begin{pmatrix} i\beta_{1}&&\\
&\ddots&\\ && i\beta_{m+2}\end{pmatrix}= \beta_{k}\, .$$ A vector
$\mu\in i\, \gt^{*}$ such that $\mu=\sum_{k=1}^{m+1} \mu_{k}\,
\widehat{x}_{k}$, in the basis $(\widehat {x}_k\equiv i\,
x_k)_{k=1,\ldots ,m+1}$, is denoted by $$ \mu= (\mu_{1}
,\mu_{2},\ldots ,\mu_{m+1})\, .$$ The restriction to $\gt$ of the
Killing form $B$ of $G$ is given by $$\forall X\in\gt\,,\,\forall
Y \in\gt\,,\quad B(X,Y)= 2\, (m+2)\,\Re\, \big(\,\mathrm{tr}
(X\,Y)\,\big)\,.$$ It is easy to verify that the scalar product on
$i \, \gt^{*}$ induced by the Killing form sign changed is given
by \begin{multline} \forall\mu=(\mu_1,\ldots ,\mu_{m+1})\in i \,
\gt^{*}\,,\,\forall \mu'=(\mu'_1,\ldots , \mu'_{m+1})\, \in i\,
\gt^*\,,\\ <\mu ,\mu'>=\frac{1}{2(m+2)}\,\sum_{k=1}^{m+1} \mu_k\,
\mu'_k -\frac{1}{2(m+2)^{2}}\,\left(\sum_{k=1}^{m+1}\mu_{k}\right)
\left(\sum_{k=1}^{m+1}\mu'_{k}\right) \,.\end{multline}
Considering the decomposition of the complexified Lie algebra of
$G$ under the action of $T$, it is easy to verify that $T$ is a
common maximal torus of $G$ and $K$, and that the respective roots
are given by
\begin{align*}
\pm(\widehat{x}_{i}-\widehat{x}_{j})\,&,\, 1\leq i<j\leq
m+1\,,&&&&
\pm\left(\widehat{x}_{i}+\sum_{k=1}^{m+1}\widehat{x}_{k}\right)\,,1\leq
i\leq m+1\,, &&\text{for } G\,,\\
\pm(\widehat{x}_{i}-\widehat{x}_{j})\,&,\, 1\leq i<j\leq m\,,&&&&
\pm\left(\widehat{x}_{m+1}+\sum_{k=1}^{m+1}\widehat{x}_{k}\right)\,,
&&\text{for } K\,.
\end{align*}
We consider as sets of positive roots
\begin{align*}\Phi_{G}^{+}&=\left\{\widehat{x}_{i}-\widehat{x}_{j}\,,\,
1\leq i\leq
m+1\,;\;\widehat{x}_{i}+\sum_{k=1}^{m+1}\widehat{x}_{k}\,,\, 1\leq
i\leq m+1\right\}\,,\\ \intertext{and} \Phi_{K}^{+}&=
\left\{\widehat{x}_{i}-\widehat{x}_{j}\,,\, 1\leq i\leq
m\,;\;\widehat{x}_{m+1}+\sum_{k=1}^{m+1}\widehat{x}_{k}\right\}
\,.
\end{align*}
Then
\begin{align*}
\delta_{G}&=\sum_{k=1}^{m+1}
(m+2-k)\,\widehat{x}_{k}=(m+1,m,\ldots,2,1)\,,\\ \intertext{and}
\delta_{K}&=\frac{1}{2}\left(\sum_{k=1}^{m}
(m+2-2k)\,\widehat{x}_{k}+2\,\widehat{x}_{m+1}\right)=
\frac{1}{2}(m,m-2,m-4\ldots,2-m,2)\,.
\end{align*}
Hence
$$\delta_{G}-\delta_{K}=\frac{1}{2}(m+2)\,\sum_{k=1}^{m}\widehat{x}_{k}
=\frac{1}{2}(m+2)(1,1,\ldots,1,0)\,,$$ so
$$\|\delta_{G}-\delta_{K}\|^{2}=\frac{m}{4}\,.$$ We now determine
the set $$\Lambda:=\{\theta\in \Phi_{G}^{+}\,;\,
<\theta,\delta_{K}>\;\;<0\}\,.$$ Recall that from the proof of
lemma~\ref{lemLambda}, if $\Lambda$ is non empty, then any
$\theta\in \Lambda$ belongs to
$\Phi_{G}^{+}\backslash\Phi_{K}^{+}$. It is then easy to verify
that the elements of $\Lambda$ are
\begin{align*}
\widehat{x}_{j}&-\widehat{x}_{m+1}\;,\; \frac{m}{2}+1\leq j\leq
m\,,
&&&&<\widehat{x}_{j}-\widehat{x}_{m+1},\delta_{K}>=\frac{1}{2(m+2)}
\left(\frac{m}{2}-j\right)\,,\\ \widehat{x}_{j}&+\sum_{k=1}^{m+1}
\widehat{x}_{k}\;,\;\frac{m}{2}+2\leq j\leq m\,,&&&& <
\widehat{x}_{j}+\sum_{k=1}^{m+1}
\widehat{x}_{k},\delta_{K}>=\frac{1}{2(m+2)}
\left(\frac{m}{2}+1-j\right)\,.
\end{align*}
So $$\sum_{\theta\in\Lambda}
<\theta,\delta_{K}>=-\frac{m^{2}}{8(m+2)}\,. $$ Hence, by
formula~(\ref{formula3}), the square of the first eigenvalue
$\lambda$ of the Dirac operator is given by $$\lambda^{2}=
\frac{m}{2}-\frac{m^{2}}{2(m+2)}+\frac{m}{2}=\frac{m+4}{m+2}\,\frac{m}{2}\,.$$

\subsection{Grassmannians $\boldsymbol{\widetilde{\mathrm{Gr}}_{4}(\mathbb{R}^{m+4})}$,
$\boldsymbol{m}$ even $\boldsymbol{\geq 4}$}

Here $G=\mathrm{Spin}_{m+4}$ and, identifying $\mathbb{R}^m$ with
the subspace of $\mathbb{R}^{m+4}$ spanned by $e_{1}, \ldots
e_{m}$, and $\mathbb{R}^{4}$ with the subspace spanned by
$e_{m+1},\ldots, e_{m+4}$, $K$ is the subgroup of $G$ defined by
$$\mathrm{Spin}_{m}\mathrm{Spin}_{4}:=\left\{\psi\in\mathrm{Spin}_{m+4}\,;\,
\psi=\varphi\phi\,,\,
\varphi\in\mathrm{Spin}_{m}\,,\,\phi\in\mathrm{Spin}_{4}
\right\}\,.$$ We consider the common torus of $G$ and $K$ defined
by $$T=\left\{ \sum_{k=1}^{\frac{m}{2}+2}\big(\cos
(\beta_{k})+\sin(\beta_{k})\,e_{2k-1}\cdot e_{2k}\big)\,;\,
\beta_{1} ,\ldots, \beta_{\frac{m}{2}+2}\, \in
\mathbb{R}\right\}\,.$$ The Lie algebra of $T$ is $$\gt=\left\{
\sum_{k=1}^{\frac{m}{2}+2}\beta_{k}\,e_{2k-1}\cdot e_{2k}\,;\,
\beta_{1} ,\ldots, \beta_{\frac{m}{2}+2}\, \in
\mathbb{R}\right\}\,.$$ We denote by $(x_1, \ldots ,
x_{\frac{m}{2}+2})$ the basis of  ${\gt}^{*}$ given by $$
x_{k}\cdot \sum_{j=1}^{\frac{m}{2}+2}\beta_{j}\,e_{2j-1}\cdot
e_{2j}=\beta_{k}\,.$$ We introduce the basis $(\widehat{x}_1,
\ldots , \widehat{x}_{\frac{m}{2}+2})$ of $i\,\gt^{*}$ defined by
$$\widehat{x}_k:=2i\,x_{k}\,,\quad k=1,\ldots,\frac{m}{2}+2\,.$$ A
vector $\mu\in i\, \gt^{*}$ such that
$\mu=\sum_{k=1}^{\frac{m}{2}+2} \mu_{k}\, \widehat{x}_{k}$, is
denoted by $$ \mu= (\mu_{1} ,\mu_{2},\ldots
,\mu_{\frac{m}{2}+2})\, .$$ The restriction to $\gt$ of the
Killing form $B$ of $G$ is given by $$ B(e_{2k-1}\cdot
e_{2k},e_{2l-1}\cdot e_{2l})=-8\,(m+2)\,\delta_{kl}\,.$$ It is
easy to verify that the scalar product on $i \, \gt^{*}$ induced
by the Killing form sign changed is given by
\begin{equation}\begin{split} \forall\mu=(\mu_1,\ldots
,\mu_{\frac{m}{2}+2})\in i \, \gt^{*}\,&,\,\forall
\mu'=(\mu'_1,\ldots , \mu'_{\frac{m}{2}+2})\, \in i\, \gt^*\,,\\
<\mu ,\mu'>&=\frac{1}{2(m+2)}\,\sum_{k=1}^{\frac{m}{2}+2} \mu_k\,
\mu'_k \,.\end{split}\end{equation} Considering the decomposition
of the complexified Lie algebra of $G$ under the action of $T$, it
is easy to verify that $T$ is a common maximal torus of $G$ and
$K$, and that the respective roots are given by
\begin{align*}
&\pm(\widehat{x}_{i}+\widehat{x}_{j})\,,\;
\pm(\widehat{x}_{i}-\widehat{x}_{j})\,,\; 1\leq i<j\leq
\frac{m}{2}+2\,,&&&\text{for } G\,,\\
&\begin{cases}\pm(\widehat{x}_{i}+\widehat{x}_{j})\,,\,
\pm(\widehat{x}_{i}-\widehat{x}_{j})\,,\,1\leq i<j\leq
\frac{m}{2}\\
\pm(\widehat{x}_{\frac{m}{2}+1}+\widehat{x}_{\frac{m}{2}+2})\,,\,
\pm(\widehat{x}_{\frac{m}{2}+1}-\widehat{x}_{\frac{m}{2}+2})\,,
\end{cases}
 &&&\text{for } K\,.
\end{align*}
We consider as sets of positive roots
\begin{align*}\Phi_{G}^{+}&=\left\{\widehat{x}_{i}+\widehat{x}_{j}\,,\,
\widehat{x}_{i}-\widehat{x}_{j}\,,\;1\leq i<j\leq
\frac{m}{2}+2\right\}\,,\\ \intertext{and} \Phi_{K}^{+}&=
\left\{\widehat{x}_{i}+\widehat{x}_{j}\,,\,
\widehat{x}_{i}-\widehat{x}_{j}\,,\;1\leq i<j\leq \frac{m}{2}\;,\;
\widehat{x}_{\frac{m}{2}+1}+\widehat{x}_{\frac{m}{2}+2}\,,\,
\widehat{x}_{\frac{m}{2}+1}-\widehat{x}_{\frac{m}{2}+2} \right\}
\,.
\end{align*}
Then
\begin{align*}
\delta_{G}&=\sum_{k=1}^{\frac{m}{2}+2}
(\frac{m}{2}+2-k)\,\widehat{x}_{k}=(\frac{m}{2}+1,\frac{m}{2},\ldots,1,0)\,,\\
\intertext{and} \delta_{K}&=\sum_{k=1}^{\frac{m}{2}}
(\frac{m}{2}-k)\,\widehat{x}_{k}+\widehat{x}_{\frac{m}{2}+1}
=(\frac{m}{2}-1,\frac{m}{2}-2,\ldots,1,0) \,.
\end{align*}
Hence
$$\delta_{G}-\delta_{K}=2\,\sum_{k=1}^{\frac{m}{2}}\widehat{x}_{k}
=2\,(1,1,\ldots,1,0,0)\,, $$ so
$$\|\delta_{G}-\delta_{K}\|^{2}=\frac{m}{m+2}\,.$$ On the other
hand, it is easy to verify that the set $$\Lambda:=\{\theta\in
\Phi_{G}^{+}\,;\, <\theta,\delta_{K}>\;\;<0\}\,,$$ has only one
element, namely
$$\widehat{x}_{\frac{m}{2}}-\widehat{x}_{\frac{m}{2}+1}\,,\text{
with }\; <\widehat{x}_{\frac{m}{2}}-\widehat{x}_{\frac{m}{2}+1},
\delta_{K}>=-1\,.$$ Hence, by formula~(\ref{formula3}), the square
of the first eigenvalue $\lambda$ of the Dirac operator is given
by $$\lambda^{2}=\frac{2m}{m+2}-\frac{2}{m+2}+\frac{m}{2}=
\frac{m^{2}+6m-4}{2(m+2)}\,.$$

\subsection{The four exceptional cases}
Note first that since all the groups $G$ we consider are simple,
their roots system are irreducible so, up to a constant, there is
only one $W_{G}$-invariant scalar product on the subspace
generated by the set of roots, cf. for instance Remark~(5.10),
\S~V in \cite{BtD}.

\noindent We use the description of root systems given in
\cite{BMP}. Those root systems are expressed in the simple root
basis $(\alpha_{i})$. Note that the $W_{G}$-invariant scalar
product $(\,,\,)$ used there is such that $(\alpha,\alpha)=2$ for
any long root $\alpha$. In order to compare it with the scalar
product $<\,,\,>$ induced by the Killing form sign-changed, we use
the ``{strange formula}'' of Freudenthal and de Vries, (cf. 47-11
in \cite{FdV}):
\begin{equation}
<\delta_{G},\delta_{G}>=\frac{1}{24}\,\dim\,G\,.
\end{equation}

\noindent To determine the set of $K$-positive roots, we use
theorem~13, theorem~14 and the proof of theorem~18 in \cite{CG}.
By those results, the set $\Phi_{K}^{+}$ may be defined as
follows. Let $\theta=\sum m_{i}\, \alpha_{i} $ be the highest
root. In all cases considered, there exists an index $j$ such that
$m_{j}=2$. Then $$\Phi_{K}^{+}=\left\{ \sum n_{i}\alpha_{i}\,;\,
n_{j}\neq 1\right\}\,.$$

\subsubsection{The symmetric space $
\mathrm{G}_{2}/\mathrm{SO}_{4}$} Using the results of pages 18 and
64 in \\\cite{BMP}, we get
$$\delta_{G}=3\,\alpha_{1}+5\,\alpha_{2}\,.$$ By the expression of
the Cartan matrix, the scalar product matrix is, in the basis
$(\alpha_{1},\alpha_{2})$, $\begin{pmatrix}
2&-1\\-1&2/3\end{pmatrix}$, hence
$$\|\delta_{G}\|_{(\,,\,)}^{2}=\frac{14}{3}\,.$$ On the other
hand, by the formula of Freudenthal and de Vries,
$$\|\delta_{G}\|_{<\,,\,>}^{2}=\frac{7}{12}\,,$$ so
$$<\;,\;>=\frac{1}{8}\,(\;,\;)\,.$$ The set of $K$-positive roots
is $$\Phi_{K}^{+}=\{2\,\alpha_{1}+3\,\alpha_{2},\alpha_{2}\}\,,$$
hence $$\delta_{K}=\alpha_{1}+2\,\alpha_{2}\,,$$ so $$\delta_{G}-
\delta_{K}=2\,\alpha_{1}+3\,\alpha_{2}\,.$$ Hence $$
\|\delta_{G}-\delta_{K}\|_{<\,,\,>}^{2}=\frac{1}{8}\,
\|\delta_{G}-\delta_{K}\|_{(\,,\,)}^{2}=\frac{1}{4}\,.$$ Finally,
it is easy to verify that the set $$\Lambda:=\{\theta\in
\Phi_{G}^{+}\,;\, <\theta,\delta_{K}>\;\;<0\}\,,$$ is empty, hence
by formula~(\ref{formula3}), the square of the first eigenvalue
$\lambda$ of the Dirac operator is given by
$$\lambda^{2}=\frac{1}{2}+1=\frac{3}{2} \,.$$

\subsubsection{The symmetric space $\mathrm{E}_{6}
/(\mathrm{SU}_{6}\mathrm{SU}_{2})$} Using the results of pages 14
and 60 in \cite{BMP}, we get
$$\delta_{G}=8\,\alpha_{1}+15\,\alpha_{2} +21\, \alpha_{3}
+15\,\alpha_{4}+8\,\alpha_{5}+11\,\alpha_{6}\,.$$ Since all roots
have same length equal to $2$, we may introduce the fundamental
weight basis $(\omega_{i})$ because
$$(\omega_{i},\alpha_{j})=\delta_{ij}\,.$$ Since $\delta_{G}=\sum
\omega_{i}$, we get $$\|\delta_{G}\|_{(\,,\,)}^{2}=78\,,$$ whereas
 by the formula of Freudenthal and de Vries,
$$\|\delta_{G}\|_{<\,,\,>}^{2}=\frac{78}{24}\,,$$ so
$$<\;,\;>=\frac{1}{24}\,(\;,\;)\,.$$ The set of $K$-positive roots
may be defined by $$\Phi_{K}^{+}=\left\{ \sum_{i=1}^{6} n_{i}\,
\alpha_{i}\,;\, n_{6}\neq 1\right\}\,.$$ Then
\begin{align*}\delta_{K}&= 3\,\alpha_{1}+5\,\alpha_{2} +6\, \alpha_{3}
+5\,\alpha_{4}+3\,\alpha_{5}+\alpha_{6}\\ &=\omega_{1}+\omega_{2}
+\omega_{3}+\omega_{4}+\omega_{5}-4\,\omega_{6}\,.
\end{align*}
Hence $$\delta_{G}-\delta_{K}= 5\,\alpha_{1}+10\,\alpha_{2} +15\,
\alpha_{3}+10\,\alpha_{4}+5\,\alpha_{5}+10\,\alpha_{6}=5\,\omega_{6}\,.$$
So $$ \|\delta_{G}-\delta_{K}\|_{<\,,\,>}^{2}=\frac{1}{24}\,
\|\delta_{G}-\delta_{K}\|_{(\,,\,)}^{2}=\frac{25}{12}\,.$$ On the
other hand it is easy to verify that the set
$$\Lambda:=\{\theta\in \Phi_{G}^{+}\,;\,
<\theta,\delta_{K}>\;\;<0\}\,,$$ has $7$ elements and that $$
\sum_{\theta\in\Lambda} <\theta,\delta_{K}> =\frac{1}{24}
\sum_{\theta\in\Lambda} (\theta,\delta_{K})=-\frac{7}{12}\,.$$ So
by formula~(\ref{formula3}), the square of the first eigenvalue
$\lambda$ of the Dirac operator is given by
$$\lambda^{2}=\frac{50}{12}-\frac{28}{12}+5=\frac{41}{6} \,.$$

\subsubsection{The symmetric space
$\mathrm{E}_{7}/(\mathrm{Spin}_{12}\mathrm{SU}_{2})$} By the
results of pages 15 and 61 in \cite{BMP}, we get
$$\delta_{G}=\frac{1}{2}\,(34\,\alpha_{1}+66\,\alpha_{2} +96\,
\alpha_{3} +75\,\alpha_{4}+52\,\alpha_{5}+27\,\alpha_{6}+
49\,\alpha_{7})\,.$$ Here again, since all roots have same length
equal to $2$, we may consider the fundamental weight basis
$(\omega_{i})$. We get
$$\|\delta_{G}\|_{(\,,\,)}^{2}=\frac{399}{2}\,,$$ whereas
 by the formula of Freudenthal and de Vries,
$$\|\delta_{G}\|_{<\,,\,>}^{2}=\frac{133}{24}\,,$$ so
$$<\;,\;>=\frac{1}{36}\,(\;,\;)\,.$$ The set of $K$-positive roots
may be defined by $$\Phi_{K}^{+}=\left\{ \sum_{i=1}^{7} n_{i}\,
\alpha_{i}\,;\, n_{1}\neq 1\right\}\,.$$ Then
\begin{align*}\delta_{K}&=\frac{1}{2}\,( 2\,\alpha_{1}+18\,\alpha_{2}
 +32\, \alpha_{3}
+27\,\alpha_{4}+20\,\alpha_{5}+11\,\alpha_{6}+17\,\alpha_{7})\\
&=-7\,\omega_{1}+\omega_{2}
+\omega_{3}+\omega_{4}+\omega_{5}+\omega_{6}+\omega_{7}\,.
\end{align*}
Hence $$\delta_{G}-\delta_{K}= 16\,\alpha_{1}+24\,\alpha_{2} +32\,
\alpha_{3}+24\,\alpha_{4}+16\,\alpha_{5}+8\,\alpha_{6}
+16\,\alpha_{7}=8\,\omega_{6}\,.$$ So $$
\|\delta_{G}-\delta_{K}\|_{<\,,\,>}^{2}=\frac{1}{36}\,
\|\delta_{G}-\delta_{K}\|_{(\,,\,)}^{2}=\frac{32}{9}\,.$$ On the
other hand it can be verified that the set $$\Lambda:=\{\theta\in
\Phi_{G}^{+}\,;\, <\theta,\delta_{K}>\;\;<0\}\,,$$ has $13$
elements and that $$ \sum_{\theta\in\Lambda} <\theta,\delta_{K}>
=\frac{1}{36} \sum_{\theta\in\Lambda}
(\theta,\delta_{K})=-\frac{41}{36}\,.$$ So by
formula~(\ref{formula3}), the square of the first eigenvalue
$\lambda$ of the Dirac operator is given by
$$\lambda^{2}=\frac{64}{9}-\frac{41}{9}+8=\frac{95}{9} \,.$$

\subsubsection{The symmetric space
$\mathrm{E}_{8}/(\mathrm{E}_{7}\mathrm{SU}_{2})$} By the results
of pages 16, 62 and 63 in \cite{BMP}, we get
$$\delta_{G}=29\,\alpha_{1}+57\,\alpha_{2} +84\, \alpha_{3}
+110\,\alpha_{4}+135\,\alpha_{5}+91\,\alpha_{6}+
46\,\alpha_{7}+68\,\alpha_{8}\,.$$ Here again, since all roots
have same length equal to $2$, we may consider the fundamental
weight basis $(\omega_{i})$. We get
$$\|\delta_{G}\|_{(\,,\,)}^{2}=620\,,$$ whereas
 by the formula of Freudenthal and de Vries,
$$\|\delta_{G}\|_{<\,,\,>}^{2}=\frac{248}{24}=\frac{31}{3}\,,$$ so
$$<\;,\;>=\frac{1}{60}\,(\;,\;)\,.$$ The set of $K$-positive roots
may be defined by $$\Phi_{K}^{+}=\left\{ \sum_{i=1}^{8} n_{i}\,
\alpha_{i}\,;\, n_{1}\neq 1\right\}\,.$$ Then
\begin{align*}\delta_{K}&=\alpha_{1}+15\,\alpha_{2}
 +28\, \alpha_{3}
+40\,\alpha_{4}+51\,\alpha_{5}+35\,\alpha_{6}+18\,\alpha_{7}
+26\,\alpha_{8}\\ &=-13\,\omega_{1}+\omega_{2}
+\omega_{3}+\omega_{4}+\omega_{5}+\omega_{6}+\omega_{7}
+\omega_{8}\,.
\end{align*}
Hence $$\delta_{G}-\delta_{K}= 28\,\alpha_{1}+42\,\alpha_{2} +56\,
\alpha_{3}+70\,\alpha_{4}+84\,\alpha_{5}+56\,\alpha_{6}
+28\,\alpha_{7}+42\,\alpha_{8}=14\,\omega_{6}\,.$$ So $$
\|\delta_{G}-\delta_{K}\|_{<\,,\,>}^{2}=\frac{1}{60}\,
\|\delta_{G}-\delta_{K}\|_{(\,,\,)}^{2}=\frac{98}{15}\,.$$ On the
other hand it can be verified that the set $$\Lambda:=\{\theta\in
\Phi_{G}^{+}\,;\, <\theta,\delta_{K}>\;\;<0\}\,,$$ has $25$
elements and that $$ \sum_{\theta\in\Lambda} <\theta,\delta_{K}>
=\frac{1}{60} \sum_{\theta\in\Lambda}
(\theta,\delta_{K})=-\frac{137}{60}\,.$$ So by
formula~(\ref{formula3}), the square of the first eigenvalue
$\lambda$ of the Dirac operator is given by
$$\lambda^{2}=\frac{196}{15}-\frac{137}{15}+14=\frac{269}{15}\,.$$

\end{document}